\documentclass{article}
\usepackage{mathrsfs}
\usepackage{amsmath}
\usepackage{amssymb, amsmath}
\usepackage{graphicx}
\catcode`\@=11 \@addtoreset{equation}{section}

\catcode`\@=12
\usepackage{color}

\newtheorem{Theorem}{Theorem}[section]
\newtheorem{Proposition}{Proposition}[section]
\newtheorem{Lemma}{Lemma}[section]
\newtheorem{Corollary}{Corollary}[section]
\newtheorem{Definition}{Definition}[section]

\newtheorem{Remark}{Remark}[section]

\newcommand{\newcom}{\newcommand}
\newcommand{\bTheorem}[1]{
\begin{Theorem} \label{T#1} }
\newcommand{\eT}{\end{Theorem}}

\newcommand{\bProposition}[1]{
\begin{Proposition} \label{P#1}}
\newcommand{\eP}{\end{Proposition}}

\newcommand{\bLemma}[1]{
\begin{Lemma} \label{L#1} }
\newcommand{\eL}{\end{Lemma}}

\newcommand{\bCorollary}[1]{
\begin{Corollary} \label{C#1} }
\newcommand{\eC}{\end{Corollary}}

\newcommand{\beq}{\begin{equation}}
\newcommand{\eeq}{\end{equation}}
\newcom{\ben}{\begin{eqnarray}}
\newcom{\een}{\end{eqnarray}}
\newcom{\beno}{\begin{eqnarray*}}
\newcom{\eeno}{\end{eqnarray*}}
\newcom{\bali}{\begin{aligned}}
\newcom{\eali}{\end{aligned}}

\newcommand{\bFormula}[1]{
\begin{equation} \label{#1}}
\newcommand{\eF}{\end{equation}}

\newcommand{\f}{\frac}
\newcommand{\Om}{\Omega}

\newcommand{\p}{\partial}

\newcommand{\vr}{\varrho}

\newcommand{\vu}{\vc{u}}

\newcommand{\vv}{\vc{v}}

\newcommand{\vw}{\vc{w}}

\newcommand{\vc}[1]{{\boldsymbol #1}}
\newcommand{\Div}{{\rm div}}
\newcommand{\Grad}{\nabla}

\newcommand{\dx}{\,{\rm d} x}

\newcommand{\dt}{\,{\rm d} t }

\newcommand{\dxdt}{\dx\dt}

\newcommand{\intTO}[1]{\int_0^T\!\!\! \int_{\Omega} #1 \ \ \dxdt}

\newcommand{\ep}{\varepsilon}

\font\F=msbm10 scaled 1000
\newcommand{\R}{\mbox{\F R}}

\newcommand{\lr}[1]{\left( #1 \right)}
\newcommand{\eq}[1]{\begin{equation}
\begin{split}
#1
\end{split}
\end{equation}}
\newcommand{\eqh}[1]{\begin{equation*}
\begin{split}
#1
\end{split}
\end{equation*}}

\newcommand\Cbox[2]{%
    \newbox\contentbox%
    \newbox\bkgdbox%
    \setbox\contentbox\hbox to \hsize{%
        \vtop{
            \kern\columnsep
            \hbox to \hsize{%
                \kern\columnsep%
                \advance\hsize by -2\columnsep%
                \setlength{\textwidth}{\hsize}%
                \vbox{
                    \parskip=\baselineskip
                    \parindent=0bp
                    #2
                }%
                \kern\columnsep%
            }%
            \kern\columnsep%
        }%
    }%
    \setbox\bkgdbox\vbox{
        \color{#1}
        \hrule width  \wd\contentbox %
               height \ht\contentbox %
               depth  \dp\contentbox
        \color{black}
    }%
    \wd\bkgdbox=0bp%
    \vbox{\hbox to \hsize{\box\bkgdbox\box\contentbox}}%
    \vskip\baselineskip%
}

\begin{document}


\title{\bf On weak solutions to the compressible inviscid two-fluid model}

\author{
Yang Li \\ School of Mathematical Sciences, \\ Anhui University, Hefei, 230601, People's Republic of China \\ lynjum@163.com\\ \\
Ewelina Zatorska \\ Department of Mathematics, \\ University College London, Gower Street, London WC1E 6BT,  United Kingdom \\
e.zatorska@ucl.ac.uk
}

\maketitle

{\centerline {\bf Abstract }}
{In three space dimensions, we consider the compressible inviscid model describing the time evolution of two fluids sharing the same velocity and enjoying the algebraic pressure closure. By employing the technique of convex integration, we prove the existence of infinitely many global-in-time weak solutions for any smooth initial data. We also show that for any piecewise constant initial densities, there exists suitable initial velocity such that the problem admits infinitely many global-in-time weak solutions that conserve the total energy. The structure of the two-fluid model is crucial in our analysis. Adaptations of the main results to other two-fluid models, like the liquid-gas flow, are available. Local-in-time existence and uniqueness of classical solutions will also be shown. 
}



\section{Introduction}

In this paper, we are concerned with the global-in-time solvability to the compressible two-fluid model system:
\begin{equation*}
\left\{\begin{aligned}
& \p_t (\alpha_{\pm} \vr_{\pm})+ \Div_x(\alpha_{\pm} \vr_{\pm} \vu)=0,\\
& \p_t [ (\alpha_{+} \vr_{+}+\alpha_{-}\vr_{-})\vu   ]+ \Div_x[(\alpha_{+} \vr_{+}+\alpha_{-}\vr_{-})\vu\otimes\vu]+\Grad_x p=\mathbf{0}, \\
& \alpha_{+}+\alpha_{-}=1,\,\,\,\alpha_{\pm}\geq 0, \\
& p=p_{+}=p_{-}. \\
\end{aligned}\right.
\end{equation*}

The model describes the motion of two immiscible compressible fluids  sharing the same velocity field and obeying the algebraic pressure closure. In the system above, $\alpha_+(t,x),\ \alpha_-(t,x)$ denote volumetric rates of presence of fluid $+$ and $-$, and so, by definition, they sum up to 1; $t\in\R^+$ and $x\in \R^3$ denote the temporal and spatial variables respectively; $\vr_+(t,x),\ \vr_-(t,x)$ denote mass densities of the constituents, $\vu(t,x)$ is the common velocity, and $p$ is the internal pressure that will be specified later on.
Throughout this paper, we neglect the effect of kinematic viscosity and immediately switch to more convenient notation using the conservative quantities:
\beq
R=\alpha_+\vr_+,\quad Q=\alpha_-\vr_-,\quad Z =\vr_+,
\eeq
in consistency with the notation from \cite{BMZ}.
The time evolution of the new ``densities" $R$ and $Q$, and the velocity field $\vu$ is now described through the following partial differential equations:
\begin{equation}\label{intr1}
\left\{\begin{aligned}
& \p_t R+\Div_x (R \vu)=0,\\
& \p_t Q+\Div_x (Q \vu)=0,\\
& \p_t [ (R+Q)\vu] +\Div_x [ (R+Q)\vu \otimes \vu ]+ \Grad_x p=\mathbf{0},\\
\end{aligned}\right.
\end{equation}
By algebraic pressure closure we mean that pressure in phase $+$ and $-$ are equal. If both of the fluids are barotropic, then 
$$p_+=\lr{\vr_+}^{\gamma_+},\quad p_-=\lr{\vr_-}^{\gamma_-}$$
with some adiabatic constants $\gamma_\pm>1$. Equality of pressures implies that $p$ can be expressed as a function of single variable $Z$
\eq{\label{intrp}
p(Z)=Z^{\gamma_+}, 
}
where $Z$ is an implicit function of $R$ and $Q$ interrelated through
\begin{equation}\label{intr2}
\left\{\begin{aligned}
& Q=\left(1-\f{R}{Z}\right) Z^{\gamma}, \quad\gamma=\gamma_{+}/\gamma_{-}, \\
& R \leq Z.
\end{aligned}\right.
\end{equation}
It was shown, see \cite{BMZ}, that $Z$ is uniquely determined by (\ref{intr2}). Thus, we may view $Z$ as a function $Z=Z(R,Q)$.

Mathematical results on the two-fluid model (\ref{intr1})-(\ref{intr2}) are mainly concentrated on its \emph{viscous} regime. In three-dimensional space, Bresch et al. \cite{BMZ} proved the existence of finite energy weak solutions to the viscous Stokes system in semi-stationary case. We refer to Novotn\'{y} and Pokorn\'{y} \cite{NP} for the existence of weak solutions to the viscous compressible two-fluid model with general pressure laws. See also \cite{LSE1} on large time behavior of weak solutions to viscous version of (\ref{intr1})-(\ref{intr2}) in one-dimensional space case. We will recall more results on other two-fluid models in Section \ref{remo}.

In this paper, we consider the \emph{global-in-time} solvability of the inviscid two-fluid model (\ref{intr1})-(\ref{intr2}) in the framework of weak solutions with large initial data. By employing the celebrated results on non-uniqueness of weak solutions to the incompressible Euler system, due to De Lellis and Sz\'ekelyhidi \cite{DS1,DS2} (see also Chiodaroli \cite{C1}, Feireisl \cite{F2}), we first prove the existence of infinitely many weak solutions to (\ref{intr1})-(\ref{intr2}) with smooth initial data. However, these weak solutions admit the initial energy jump on account of the convex integration scheme. Inspired by Feireisl et al. \cite{FKKM} (see also Luo et al. \cite{LXX}), we then show the existence of infinitely many weak solutions which conserve the total energy for piecewise constant initial densities and suitably chosen initial velocity. It turns out that the arguments used for (\ref{intr1})-(\ref{intr2}) can be proceeded in the same way for some other two-fluid models, see Section \ref{remo}. In Section \ref{local}, we show that the Cauchy problem for (\ref{intr1})-(\ref{intr2}) admits \emph{local-in-time} existence and uniqueness of classical solutions. The main idea is to symmetrize the system in such a way that it enjoys similar structure to compressible Euler system, as observed by Ruan and Trankhinin \cite{RuTr}. This gives the local-in-time well-posedness in the framework of classical solutions.   

We remark that the verification of non-uniqueness results for the incompressible Euler system dates back to the papers of Scheffer \cite{Sche} and Shnirelman \cite{Shir} in two space dimensions. These results were greatly improved by De Lellis and Sz\'ekelyhidi \cite{DS1,DS2} by the technique of convex integration. Later, the convex integration scheme was adapted to the compressible Euler system by Chiodaroli \cite{C1} and Feireisl \cite{F2} of an abstract Euler-type system. For more results on its applications on models arising from fluid dynamics, we refer to \cite{CM1,CFK,CFG,DFM,F1,FGS,FL1,Shvy}, among others.

\section{The main results}

To fix ideas, we assume the fluids occupy a bounded Lipschitz domain $\Om\subset \R^3$. System (\ref{intr1}) is supplemented with the initial conditions:
\beq\label{intr3}
[R,Q,\vu](0,\cdot)=[R_0,Q_0,\vu_0](\cdot) \,\,\text{ in }\Om,
\eeq
and the impermeability boundary conditions:
\beq\label{intr4}
\vu\cdot \vc{n}=0 \,\, \text{ on }\p \Om,
\eeq
where $\vc{n}$ is the unit outward normal to $\p \Om$.

We are now ready to give the definition of weak solutions as follows.
\begin{Definition}\label{def:1}
Let $T\in (0,\infty)$. A triple $[R,Q,\vu]$ is said to be a weak solution to (\ref{intr1})-(\ref{intr2}),  with the initial conditions  (\ref{intr3}) and the boundary conditions (\ref{intr4}) in $(0,T)\times \Om$ provided that
\begin{itemize}
\item {\[
R(t,x),\,Q(t,x)\geq 0
\]\\
  for a.e. $(t,x)\in (0,T)\times \Om$;}

\item { \[
\intTO{\left(R \p_t \phi +R \vu\cdot \Grad_x \phi\right)} +\int_{\Om} R_0 \phi(0,\cdot)\dx=0
\]\\
      for any $\phi\in C_c^{\infty}( [0,T)\times \overline{\Om})$;}

\item { \[
\intTO{\left(Q \p_t \phi +Q \vu\cdot \Grad_x \phi\right)} +\int_{\Om} Q_0 \phi(0,\cdot)\dx=0
\]\\
      for any $\phi\in C_c^{\infty}( [0,T)\times \overline{\Om})$;}

\eqh{
\intTO{\Big( (R+Q) \vu\cdot\p_t \vc{\phi} +(R+Q) \vu\otimes \vu: \Grad_x \vc{\phi}+ Z^{\gamma_{+}}\Div_x \vc{\phi}\Big)}\\
    +\int_{\Om} (R_0+Q_0) \vu_0\cdot \vc{\phi}(0,\cdot)\dx=0
}
      for any $\vc{\phi}\in C_c^{\infty}( [0,T)\times \overline{\Om};\R^3),\,\,\vc{\phi} \cdot \vc{n}|_{\p \Om}=0$.
\end{itemize}
\end{Definition}

Our first result is concerned with the existence of global-in-time weak solutions with general initial data.
\begin{Theorem}\label{TH1}
Let $\Om\subset \R^3$ be a bounded Lipschitz domain and $\gamma_{\pm}>1$. Suppose that
\[
R_0>\underline{R}>0,\,\,Q_0>\underline{Q}>0,
\]
\[
R_0,Q_0 \in C^3(\overline{\Om}),\, \vu_0 \in C^3(\overline{\Om};\R^3),\,\,\vu_0\cdot \vc{n}|
_{\p \Om}=0.
\]

Then there exist infinitely many weak solutions to (\ref{intr1})-(\ref{intr2}) in the sense of Definition \ref{def:1}, for any $T\in (0,\infty)$.
\end{Theorem}

By setting $\alpha:= \f{R}{Z}$, one formally derives from (\ref{intr1}) the energy identity (see Lemma 2.2 in \cite{BMZ}):
\beq\label{lu1}
\f{d}{dt} \int_{\Om}
\left[
\f{1}{2}(R+Q) |\vu|^2   + \f{1}{ \gamma_{+}-1 } \left(\f{R}{\alpha} \right)^{ \gamma_{+} }  \alpha
+ \f{1}{ \gamma_{-}-1 } \left(\f{Q}{1-\alpha} \right)^{ \gamma_{-} }  (  1-\alpha  )
\right]
\dx=0.
\eeq
We remark that the weak solutions obtained in Theorem \ref{TH1} admit the drawback of initial energy jump due to the method of convex integration. Thus, the energy inequality is not satisfied. In the next theorem, suitable initial data are chosen to avoid this.

\begin{Theorem}\label{TH2}
Let $\Om$ and $\gamma_{\pm}$ be as in Theorem \ref{TH1}. Assume that the initial densities $R_0$ and $Q_0$ are piecewise constant and bounded, i.e., there exist at most countably open sets $\{\Om_i \}$ such that
\[
\Om =\cup_{i}\Om_i,\,\,\Om_i \subset \Om,\,\,\Om_i \cap \Om_j=\emptyset \text{ if }i\neq j,\,\,|\p \Om_i|=0,
\]
\[
R_0|_{\Om_i}=R_0^i ,\,\,Q_0|_{\Om_i}=Q_0^i,
\]
\[
0< \inf_{i}R_0^i \leq \sup_{i}R_0^i <\infty, \,\, 0< \inf_{i}Q_0^i \leq \sup_{i}Q_0^i <\infty.
\]

Then there exists $\vu_0\in L^{\infty}(\Om;\R^3)$ such that the problem (\ref{intr1})-(\ref{intr2}) admits infinitely many weak solutions in the sense of Definition \ref{def:1} emanating from the same initial data $[R_0,Q_0,\vu_0]$. Furthermore, these solutions comply with the conservation of total energy \eqref{lu1}.
\end{Theorem}
\begin{Remark}
With obvious modifications in the proof, the conclusions of Theorems \ref{TH1}-\ref{TH2} hold for the $N$-dimensional space with $N\geq 2$. Here, we choose the physically relevant three-dimensional space for definiteness.
\end{Remark}

The rest of this paper is arranged as follows. In Section \ref{mawe}, we show the existence of infinitely many weak solutions with general initial data. By choosing piecewise constant initial densities and suitable initial velocity, we prove the existence of infinitely many weak solutions satisfying the energy equality, see Section \ref{maad}. Further discussions will be presented in Section \ref{fur} concerning the adaptations to other two-fluid models and local-in-time existence and uniqueness of classical solutions.

\section{Infinitely many weak solutions}\label{mawe}

The present section is dedicated to the proof of Theorem \ref{TH1}. The arguments used here are adaptations of the non-uniqueness results for compressible Euler system from \cite{CFK}. 

Our purpose first is to rewrite system \eqref{intr1} as a variant of incompressible Euler system, to be able to apply the convex integration machinery. We employ the Helmholz decomposition and we identify two potentials $\Psi_1$, $\Psi_2$ and two solenoidal vector fields $\vv_1 $, $\vv_2 $ such that
\beq\label{Inf1}
R \vu =\vv_1 +\Grad_x \Psi_1,\,\,\Div_x \vv_1 =0,
\eeq
\beq\label{Inf2}
Q \vu =\vv_2 +\Grad_x \Psi_2,\,\,\Div_x \vv_2 =0.
\eeq
To proceed, we use the following ansatz:
\begin{equation}\label{Inf3}
\left\{\begin{aligned}
& R,\,Q \in C^2([0,T]\times \overline{\Om}),\,\, R(0,\cdot)=R_0(\cdot),\,\,Q(0,\cdot)=Q_0(\cdot),\\
& R(t,x),\,Q(t,x) \geq \underline{C}>0,  \text{ for any } (t,x)\in [0,T]\times \overline{\Om}.\\
\end{aligned}\right.
\end{equation}
Then, recalling the continuity equations (\ref{intr1})$_1$ and (\ref{intr1})$_2$, one finds that the potential functions $\Psi_1$ and $\Psi_2$ are unique solutions to the  following Neumann problems:
\begin{equation}\label{Inf4}
\left\{\begin{aligned}
& -\Delta \Psi_1 (t,\cdot)=\p_t R(t,\cdot),\\
& \Grad_x \Psi_1 \cdot \vc{n} |_{\p \Om}=0, \\
& \int_{\Om}\Psi_1 (t,\cdot) \dx=0,\\
\end{aligned}\right.
\end{equation}
\begin{equation}\label{Inf5}
\left\{\begin{aligned}
& -\Delta \Psi_2 (t,\cdot)=\p_t Q(t,\cdot),\\
& \Grad_x \Psi_2 \cdot \vc{n} |_{\p \Om}=0, \\
& \int_{\Om}\Psi_2 (t,\cdot) \dx=0.\\
\end{aligned}\right.
\end{equation}
It follows from (\ref{Inf1}) and (\ref{Inf2}) that
\beq\label{Inf6}
(R+Q) \vu =(\vv_1 +\vv_2) + \Grad_x (\Psi_1  + \Psi_2)=: \vv+ \Grad_x \Psi.
\eeq
Consequently, it remains to prove the existence of infinitely many weak solutions to
\begin{equation}\label{Inf7}
\left\{\begin{aligned}
& \p_t \vv+\Div_x\left(\f{(\vv+\Grad_x \Psi) \otimes (\vv+\Grad_x \Psi)} {R+ Q}\right)
 +\Grad_x \left(
 Z^{ \gamma_{+} }(R,Q) +\p_t \Psi-\f{2}{3} \Lambda
  \right) =\mathbf{0}, \\
 & \Div_x \vv=0, \\
& \vv \cdot \vc{n} | _{\p \Om}=0,\\
& \vv(0,\cdot)=\vv_0:=R_0\vu_0-\Grad_x \Psi_1(0,\cdot)+ Q_0 \vu_0 -\Grad_x \Psi_2 (0,\cdot).\\
\end{aligned}\right.
\end{equation}
In (\ref{Inf7})$_1$, $\Lambda$ is a positive constant to be determined later on. 

Notice that due to the ansatz \eqref{Inf3}, the triple $[R,Q,\vu]$ with
\[
\vu=\f{\vv + \Grad_x \Psi}{R+Q}
\]
solves the original problem (\ref{intr1}) with the initial-boundary conditions (\ref{intr3})-(\ref{intr4}) as long as $\vv$ is a solution to (\ref{Inf7}).
Notice also  that (\ref{Inf7}) is of relevant form, similar to the incompressible Euler system. The existence of infinitely many weak solutions to such a type of systems has been obtained by Feireisl \cite{F2} in a quite general setting. Inspired by \cite{CFK,F2}, we set the kinetic energy as
\eq{\label{kin}
e:=\Lambda - \f{3}{2} \Big(
 Z^{ \gamma_{+} }(R,Q) +\p_t \Psi
\Big).
}
The crucial concept in the convex integration scheme lies in the choice of the set of subsolutions (see \cite{DS2,CFK,F1,F2,FGS,FL1}). 
\begin{Definition}\label{def:sub}
Let $\R^{3\times 3}_{sym,0}$ denote the space of symmetric $3\times 3$ matrices with zero trace, and  $\lambda_{max}[\mathbb{V}]$ denote the maximal eigenvalue of $\mathbb{V}\in \R^{3\times 3}_{sym}$. 
We will call $X_0$ the set of subsolutions if
\[
X_0:=\Big\{\vv\,\,|\vv\in L^{\infty}((0,T)\times \Om;\R^3) \cap C_{weak}([0,T];L^2(\Om;\R^3))
\cap C^1((0,T)\times \Om;\R^3),
\]
such that $\vv$ satisfies 
\begin{equation*}
\left\{\begin{aligned}
& \p_t \vv+\Div_x \mathbb{U}=\mathbf{0},\\
& \Div_x \vv=0,\\
& \vv(0,\cdot)=\vv_0, \\
& \vv \cdot \vc{n}|_{\p \Om}=0, \\
\end{aligned}\right.
\end{equation*}
\centerline{for some $\mathbb{U}\in C^1((0,T)\times \Om;\R^{3\times 3}_{sym,0})\cap L^{\infty}((0,T)\times \Om;\R^{3\times 3}_{sym,0})$,}
\[
\f{3}{2}\lambda_{max}\left[\f{(\vv+\Grad_x \Psi) \otimes (\vv+\Grad_x \Psi)} {R+Q}-\mathbb{U}  \right] < e\,\, \text{ in }(0,T)\times \Om \Big\}.
\]
\end{Definition}

To make sure that $X_0$ is non-empty, we first make use of an observation from \cite{LSE1}. Suppose that
\beq\label{Inf8}
0< \underline{R}\leq R (t,x) \leq \overline{R}< \infty,\,\,
0< \underline{Q}\leq Q (t,x) \leq \overline{Q}< \infty,
\eeq
for any $(t,x)\in [0,T]\times \overline{\Om}$. Then it follows from the structural relations (\ref{intr2}) that there exist constants $\underline{Z}, \,\overline{Z}$ such that
\beq\label{Inf9}
0< \underline{Z}\leq Z (t,x) \leq \overline{Z}< \infty,\,\,
\eeq
for any $(t,x)\in [0,T]\times \overline{\Om}$. Indeed, the lower bound of $Z$ follows directly from (\ref{intr2})$_2$ and it is enough to take $\underline{Z}=\underline{R}$. Notice also that $Z$ must be bounded from above, i.e.,
\[
Z(t,x)\leq \max\left\{2\overline{R},(2\overline{Q})^{1/{\gamma}}\right\},\,\,\text{ for any } (t,x)\in [0,T]\times \overline{\Om}.
\]
Otherwise, assume that there exists $(t_0,x_0)\in [0,T]\times \overline{\Om}$ such that
\[
Z(t_0,x_0) > \max\left\{2\overline{R},(2\overline{Q})^{1/{\gamma}}\right\},
\]
then
\[
\overline{Q}\geq Q(t_0,x_0)=\left(1-\f{R(t_0,x_0)}{Z(t_0,x_0)}\right)Z^{\gamma} (t_0,x_0)  \geq \f{1}{2}Z^{\gamma}  (t_0,x_0) >\overline{Q},
\]
which is a contradiction.

Now, let $E(\vv,\mathbb{U}):=\lambda_{max}(\vv \otimes \vv- \mathbb{U})$ for $\vv\in \R^3,\mathbb{U}\in \R^{3\times 3}_{sym,0}$. It follows from Lemma 3 in \cite{DS2} that the mapping $E(\cdot,\cdot):\R^3 \times \R^{3\times 3}_{sym,0}\rightarrow \R$ is convex. Thus,
\[
\lambda_{max}\left[\f{(\vv_0+\Grad_x \Psi) \otimes (\vv_0+\Grad_x \Psi)} {R+Q} \right]
\]
is bounded. Since from  (\ref{Inf3})-(\ref{Inf5}) $\p_t \Psi$ is bounded as well, and thanks to (\ref{Inf8})-(\ref{Inf9}), we may now fix $\Lambda$ sufficiently large such that 
\beq\label{Inf10}
\f{3}{2}\lambda_{max}\left[\f{(\vv_0+\Grad_x \Psi) \otimes (\vv_0+\Grad_x \Psi)} {R+Q} \right]
+
 \f{3}{2} \Big(
 Z^{ \gamma_{+} }(R,Q) +\p_t \Psi
\Big)
< \Lambda
\eeq
holds in $ (0,T) \times \Om$. We thus found that $\vv_0$ together with $\mathbb{U}=0$ belongs to the space of subsolutions $X_0$.

To proceed, we recall the crucial oscillatory lemma. 
\begin{Lemma}[Lemma 3.2, \cite{CFK}]\label{osci}
Let $\tilde{\vr}$, $\mathbf{V} $ satisfy
\[
0<\underline{\vr}\leq \tilde{\vr}\leq \overline{\vr}<\infty,\,\,\tilde{\vr} \in
C^1([0,T]\times \overline{\Om}),
\]
\[
\widetilde{\mathbf{V}  } \in C^1([0,T]\times \overline{\Om};\R^3).
\]
Suppose in addition that
\[
\vv \in  C_{weak}([0,T];L^2(\Om;\R^3))
\cap C^1((0,T)\times \Om;\R^3),
\]
solves the linear system
\begin{equation*}
\left\{\begin{aligned}
& \p_t \vv+\Div_x \mathbb{W}=\mathbf{0},\\
& \Div_x \vv=0,\\
\end{aligned}\right.
\end{equation*}
in $(0,T)\times \Om$ for some $\mathbb{W}$ belonging to $C^1((0,T)\times \Om;\R^{3\times 3}_{sym,0})$ such that
\[
\f{3}{2}\lambda_{max}\left[   \f{(\vv+\widetilde{\mathbf{V} } ) \otimes (\vv+\widetilde{\mathbf{V} } )}{\tilde{\vr}}-\mathbb{W}\right]< \tilde{e}
\]
in $(0,T)\times \Om$, for some $\tilde{e}\in C([0,T]\times \overline{\Om}),\,\tilde{e}>0$. 

Then there exist two sequences
\[
\{ \vv_n \}_{n=1}^{\infty}\subset C_c^{\infty}((0,T)\times \Om;\R^3),\,\, \{ \mathbb{W}_n \}_{n= 1}^{\infty} \subset
C_c^{\infty}((0,T)\times \Om;\R^{3\times 3}_{sym,0})
\]
such that
\begin{equation}\label{Inf11}
\left\{\begin{aligned}
& \p_t \vv_n +\Div_x \mathbb{W}_n=0,\\
& \Div_x \vv_n=0, \\
& \f{3}{2} \lambda_{max}\left[   \f{(\vv+\widetilde{\mathbf{V} } +\vv_n) \otimes (\vv+\widetilde{\mathbf{V} } +\vv_n)}{\tilde{\vr}}-(\mathbb{W}+\mathbb{W}_n)\right]< \tilde{e}\text{ in }(0,T)\times \Om,\\
& \vv_n \rightarrow \mathbf{0} \text{ in } C_{weak}([0,T];L^2(\Om;\R^3)),\\
& \liminf_{n\rightarrow \infty}\intTO{\f{|\vv_n|^2}{\tilde{\vr}}}
\geq C_{\ast}  \intTO{ \left( \tilde{e}-\f{1}{2}\f{| \vv+\widetilde{\mathbf{V} }|^2}{\tilde{\vr}} \right)^2}.\\
\end{aligned}\right.
\end{equation}
Here, $C_{\ast} $ is a positive constant depending solely on $\tilde{\vr},\widetilde{\mathbf{V} }$ and $\tilde{e}$.
\end{Lemma}
This lemma may be viewed as a variant of Proposition 3 from \cite{DS2}  in the context of incompressible Euler equations (see also Chiodaroli \cite{C1}).

Let $\vv\in X_0$ with the associated matrix-valued function $\mathbb{U}$. Using again the elementary inequality from Lemma 3 in \cite{DS2} we get that
\beq\label{Inf12}
\f{1}{2}\f{|\vv +\Grad_x \Psi|^2}
{R+Q}\leq
\f{3}{2}\lambda_{max}\left[\f{(\vv+\Grad_x \Psi) \otimes (\vv+\Grad_x \Psi)} {R+Q}-\mathbb{U}  \right]< e.
\eeq
We then conclude from (\ref{Inf12}) that the set of subsolutions is bounded in $L^{\infty}((0,T)\times \Om;\R^3)$. As an immediate consequence, we realize that the space of subsolutions $X_0$ with respect to the topology of $C_{weak}([0,T];L^2(\Om;\R^3)$) is metrizable, see for example \cite{DS2}.
 Let $\overline{X_0}$ be the closure of $X_0$ with respect to this metrizable topology. In other words, $\overline{X_0}$ becomes a complete metric space. This will allow us to apply  Baire's category theorem in the end of this section. 

We are now in a position to complete the proof of Theorem \ref{TH1}, adapting the arguments from \cite{CFK,F2,FL1}.

\emph{Proof of Theorem \ref{TH1}.} By setting the functional
\[
I[\vv]:=\int_0^T\int_{\Om}\left(\f{1}{2}\f{|\vv+\Grad_x \Psi|^2}{R+Q}
-e \right)\dxdt,
\]
in agreement with the definition of subsolutions, we may regard $I$ as a mapping ranging from $\overline{X_0}$ to $(-\infty,0]$. Now the central task is to show that
\beq\label{Inf13}
I[\vv]=0,\text{ if }I\text{ is continuous at }\vv\in \overline{X_0}.
\eeq
Suppose that $\vv\in \overline{X_0}$ is a continuity point of $I$. Then there exist a sequence $\{\vv^{(k)} \}_{k=1}^{\infty}\subset X_0$ such that
\[
\vv^{(k)} \rightarrow \vv \text{ in } C_{weak}([0,T];L^2(\Om;\R^3))
\]
and
\[
I[\vv^{(k)}]\rightarrow I[\vv]
\]
as $k\rightarrow \infty$. Let $\{ \mathbb{U}^{(k)} \}_{k=1}^{\infty}\subset C^1((0,T)\times \Om;\R^{3\times 3}_{sym,0})\cap L^{\infty}((0,T)\times \Om;\R^{3\times 3}_{sym,0})$ be the sequence of related fluxes in light of the Definition \ref{def:sub}. Then
\[
\f{3}{2}\lambda_{max}\left[\f{(\vv^{(k)}+\Grad_x \Psi) \otimes (\vv^{(k)}+\Grad_x \Psi)} {R+Q}-\mathbb{U}^{(k)}\right] < e-\ep_k\,\, \text{ in }(0,T)\times \Om
\]
for a suitable vanishing sequence $\{ \ep_k \}_{k=1}^{\infty}$.
For a fixed $k$ we may now apply Lemma \ref{osci} with 
\[
\tilde{\vr}=R+Q,\quad\vv=\vv^{(k)} ,\quad \widetilde{\mathbf{V}  }=\Grad_x \Psi,\quad
\mathbb{W}=\mathbb{U}^{(k)},\quad \tilde{e}=e-\ep_k,
\]
we conclude that there exist two sequences
\[
\{ \vv_n^{(k)} \}_{n=1}^{\infty}\subset C_c^{\infty}((0,T)\times \Om;\R^3),\quad \{ \mathbb{W}_n ^{(k)} \}_{n= 1}^{\infty} \subset
C_c^{\infty}((0,T)\times \Om;\R^{3\times 3}_{sym,0})
\]
such that
\begin{equation}\label{Inf14}
\left\{\begin{aligned}
& \p_t \vv_n^{(k)} +\Div_x \mathbb{W}_n^{(k)}=0,\\
& \Div_x \vv_n^{(k)}=0, \\
& \f{3}{2} \lambda_{max}\left[   \f{(\vv^{(k)}+\Grad_x \Psi +\vv_n^{(k)}) \otimes (\vv^{(k)}+\Grad_x \Psi +\vv_n^{(k)})}{R+Q}-(\mathbb{U}^{(k)}+\mathbb{U}^{(k)}_n)\right]< e-\ep_k \text{ in }(0,T)\times \Om,\\
& \vv_n^{(k)} \rightarrow \mathbf{0} \text{ in } C_{weak}([0,T];L^2(\Om;\R^3)),\\
& \liminf_{n\rightarrow \infty}\int_0^T \int_{\Om}\f{|\vv_n^{(k)}|^2}{R+Q}\dxdt
\geq C_{\ast} \int_0^T\int_{\Om} \left( e-\ep_k-\f{1}{2}\f{| \vv^{(k)}+\Grad_x \Psi|^2}{R+Q} \right)^2 \dxdt.\\
\end{aligned}\right.
\end{equation}
In particular, $C_{\ast} $ is independent of $k$ and $n$. 
Upon setting
\[
\vw _n^{(k)}:=\vv^{(k)}+\vv_n^{(k)},\,\,
\mathbb{H}_n^{(k)}:=\mathbb{U}^{(k)}+\mathbb{U}^{(k)}_n,\,\,
\]
it follows that
\begin{equation}\label{Inf15}
\left\{\begin{aligned}
& \p_t \vw _n^{(k)} +\Div_x \mathbb{H}_n^{(k)}=0,\\
& \Div_x \vw _n^{(k)}=0, \\
& \f{3}{2} \lambda_{max}\left[   \f{(\vw _n^{(k)}+\Grad_x \Psi ) \otimes (\vw _n^{(k)}+\Grad_x \Psi )}{R+Q}-\mathbb{H}_n^{(k)}\right]< e-\ep_k \text{ in }(0,T)\times \Om.\\
\end{aligned}\right.
\end{equation}
As a consequence, we find that $\vw_n^{(k)}\in X_0$ for any $n,k \geq 1$.

Making use of $(\ref{Inf14})_4$ and Cauchy-Schwarz's inequality,
\eq{\label{Inf16}
&\liminf_{k\rightarrow \infty}I[\vw_k^{(k)}]
=\liminf_{k\rightarrow \infty}\intTO{
\left(\f{1}{2}\f{|\vv^{(k)}+\vv_k^{(k)}+\Grad_x \Psi|^2}{R+Q}
-e \right)}
\\
&=\lim_{k\rightarrow \infty}\int_0^T\int_{\Om}
\left(\f{1}{2}\f{|\vv^{(k)}+\Grad_x \Psi|^2}{R+Q}
-e \right)\dxdt
+\liminf_{k\rightarrow \infty} \int_0^T\int_{\Om} \f{1}{2}\f{|\vv_k^{(k)}|^2}{R+Q} \dxdt
\\
&\geq I[\vv]+\f{C_{\ast}}{2}\liminf_{k\rightarrow \infty}\int_0^T\int_{\Om}
\left( e-\ep_k-\f{1}{2}\f{| \vv^{(k)}+\Grad_x \Psi|^2}{R+Q} \right)^2
\dxdt
\\
&\geq I[\mathbf{w}]+\f{C_{\ast}}{2T}\liminf_{k\rightarrow \infty}\left\{\int_0^T\int_{\Om}
\left( e-\ep_k-\f{1}{2}\f{| \vv^{(k)}+\Grad_x \Psi|^2}{R+Q} \right)
\dxdt\right\}^2
\\
&=I[\vv]+\f{C_{\ast}}{2T}\left(  I[\vv] \right)^2.
}
Notice that $I[\vw_k^{(k)}]\rightarrow I[\vv]$ because $\vw_k^{(k)}\rightarrow \vv$ in $C_{weak}([0,T];L^2(\Om;\R^3))$ as $k\rightarrow \infty$. Then it follows from (\ref{Inf16}) that $I[\vv]=0$, which verifies (\ref{Inf13}).

Hence, for any $\vv \in \overline{X_0}$ which is a continuity point of the functional $I$, the inequality \eqref{Inf12} must be an equality. Recalling again  Lemma 3 from \cite{DS2}, this happens only if
\[
\mathbb{U}=\f{  (\vv+\Grad_x \Psi)\otimes (\vv+\Grad_x \Psi)         }{R+Q}
-\f{1}{3}     \f{ |\vv+\Grad_x \Psi|^2}{R+Q}            \mathbb{I}_3,
\]
where the associated flux $\mathbb{U}$ may be obtained as the weak limits of the corresponding fluxes in $X_0$. 
Passing to the limit and using the definition of  $e$, \eqref{kin}, we deduce that $\vv$ solves (\ref{Inf7}). 

Finally, observe that the functional $I$ is lower semi-continuous on the complete metric space $\overline{X_0}$. In accordance with Baire's category theorem, the points of continuity of $I$ admit infinite cardinality. This finishes the proof of Theorem \ref{TH1}. $\Box$

\section{Infinitely many admissible weak solutions}\label{maad}

The aim of this section is to pick up suitable initial data such that the problem (\ref{intr1})-(\ref{intr4}) admits infinitely many weak solutions which conserve the total energy \eqref{lu1}. These weak solutions are usually termed as \emph{admissible weak solutions}.

\emph{Proof of Theorem \ref{TH2}.} 

\noindent \emph{Step 1.} We begin with the special case that the initial densities 
$$R_0\equiv R=const.>0,\quad Q_0 \equiv Q=const.>0.$$ 
Let $Z(R,Q)$ be the positive constant determined by the relation (\ref{intr2}). It follows from Theorem 13.6.1 in  \cite{F2}  
that there exists $\vc{m}_0\in L^{\infty}(\Om;\R^3)$ and a positive constant $\chi$ such that the problem
\begin{equation}\label{lu4}
\left\{\begin{aligned}
& \Div_x \vc{m}=0,\\
& \p_t \vc{m}+\Div_x\left(\f{\vc{m} \otimes \vc{m}}{R+Q}-\f{1}{3}\f{|\vc{m}|^2}{R+Q}\mathbb{I}_3\right)=
\mathbf{0},\\
& \vc{m}(0,\cdot)=\vc{m}_0,\\
& \vc{m} \cdot \vc{n}|_{\p \Om}=0,\\
\end{aligned}\right.
\end{equation}
admits infinitely many weak solutions $\vc{m}$ in $(0,T)\times \Om$ obeying
\begin{equation}\label{lu5}
\left\{\begin{aligned}
&   \vc{m}\in  L^{\infty}((0,T)\times \Om;\R^3) \cap C_{weak}([0,T];L^2(\Om;\R^3)),\\
& \f{1}{2}\f{|\vc{m}|^2}{R+Q}=\chi- \f{3}{2} Z^{\gamma_{+}  } (R,Q) \,\,\text{ for a.e.  }t\in (0,T),x\in \Om,\\
& \f{1}{2}\f{|\vc{m}_0|^2}{R+Q}=\chi- \f{3}{2} Z^{\gamma_{+}  } (R,Q) .\\
\end{aligned}\right.
\end{equation}
The role of $\chi$ is similar to $\Lambda$ in (3.7).
It follows immediately from (\ref{lu5}) that
\beq\label{lu51}
\int_{\Om}  \f{1}{2}\f{|\vc{m}|^2}{R+Q}(t,x)\dx=\int_{\Om} \f{1}{2}\f{|\vc{m}_0|^2}{R+Q}\dx,
\eeq
for a.e. $t \in (0,T)$.  In addition, (\ref{lu4})$_1$ and (\ref{lu4})$_2$ are satisfied in the sense of distributions, i.e.,
\beq\label{lu6}
\intTO{\vc{m}\cdot \Grad_x \phi} =0,  \text{ for any } \phi\in C^{\infty}( [0,T]\times \overline{\Om}),
\eeq
\beq\label{lu7}
\intTO{\left[ \vc{m}\cdot \p_t \vc{\phi} +\left(\f{\vc{m} \otimes \vc{m}}{R+Q}-\f{1}{3}\f{|\vc{m}|^2}{R+Q}\mathbb{I}_3\right):\Grad_x \vc{\phi}\right]} +\int_{\Om} \vc{m}_0\cdot\vc{\phi}(0,\cdot)\dx=0,
\eeq
for any $\vc{\phi}\in C_c^{\infty}( [0,T)\times \overline{\Om};\R^3)$. It should be emphasized that (\ref{lu7}) holds for any smooth test function without requiring the normal trace of $\vc{\phi}$ to be zero, as observed by Feireisl et al. in \cite{FKKM} in the context of complete Euler system. Thus, by setting the velocity to be $\vu:= \f{\vc{m}}{R+Q}$, we conclude from  (\ref{lu6}) that
\beq\label{lu8}
\intTO{\left(R \p_t \phi +R \vu\cdot \Grad_x \phi\right)} +\int_{\Om} R \phi(0,\cdot)\dx=0,
\eeq
\beq\label{lu81}
\int_0^T \int_{\Om}\left(Q \p_t \phi +Q \vu\cdot \Grad_x \phi\right) \dxdt +\int_{\Om} Q \phi(0,\cdot)\dx=0,
\eeq
for any $\phi \in C_c^{\infty}( [0,T)\times \overline{\Om})$ and from (\ref{lu5})$_2$  and (\ref{lu7}) that
\eq{\label{lu9}
\intTO{\left\{ (R+Q) \vu\cdot \p_t \vc{\phi} +\left[ (R+Q) \vu\otimes \vu +
Z^{\gamma_{+}  } (R,Q) \mathbb{I}_3
-\f{2}{3}\chi\,\mathbb{I}_3\right]:\Grad_x \vc{\phi}\right\}}\\
+\int_{\Om} (R+Q)\vu_0 \cdot \vc{\phi}(0,\cdot)  \dx=0,
}
for any $\vc{\phi}\in C_c^{\infty}( [0,T)\times \overline{\Om};\R^3)$. Observe that the $\chi$-dependent term vanishes after restricting the test functions to be $\vc{\phi} \cdot \vc{n}|_{\p \Om}=0$. Finally, it follows from (\ref{lu51}) and the fact that $R,\,Q$ are positive constants that
\eq{\label{lu91}
 &\int_{\Om}
\left[
\f{1}{2}(R+Q) |\vu|^2   + \f{1}{ \gamma_{+}-1 } \left(\f{R}{\alpha} \right)^{ \gamma_{+} }  \alpha
+ \f{1}{ \gamma_{-}-1 } \left(\f{Q}{1-\alpha} \right)^{ \gamma_{-} }  (  1-\alpha  )
\right](t)
\dx\\
&=
 \int_{\Om}
\left[
\f{1}{2}(R_0+Q_0) |\vu_0|^2   + \f{1}{ \gamma_{+}-1 } \left(\f{R_0}{\alpha_0} \right)^{ \gamma_{+} }  \alpha
+ \f{1}{ \gamma_{-}-1 } \left(\f{Q_0}{1-\alpha_0} \right)^{ \gamma_{-} }  (  1-\alpha_0  )
\right]
\dx,
}
for a.e. $t\in (0,T)$. Recall that we defined $\alpha=\frac{R}{Z}$, so by taking $R$, $Z$ constant, we guarantee that also $\alpha=\alpha_0=conts.$

Combining (\ref{lu8})-(\ref{lu91}), we have verified that for any  constants positive densities $R_0,Q_0$ and $R,Q$ such that  $R_0\equiv R,\, Q_0 \equiv Q$, there exists $\vu_0 \in L^{\infty}(\Om;\R^3)$ such that there exist infinitely many weak solutions $[R,Q,\vu]$ to the problem (\ref{intr1})-(\ref{intr4}) in $(0,T)\times \Om$ conserving the total energy \eqref{lu91}.

\smallskip

\noindent \emph{Step 2.} We now turn to more general case when the initial densities are piecewise constant as assumed in Theorem \ref{TH2}. Employing the result of Step 1 in each $\Om_i$, we obtain a suitable initial momentum $\vc{m}_0^i\in L^{\infty}(\Om_i;\R^3)$ such that the problem (\ref{intr1}) admits infinitely many weak solutions $[R_0^i,Q_0^i,\vc{m}^i]$ in $(0,T)\times \Om_i$ with initial data $[R_0^i,Q_0^i,\vu_0^i]$. Here, $\vu_0^i:= \f{\vc{m}_0^i}{R_0^i+Q_0^i}$. Consequently, by setting $$\vc{m}_0|_{\Om_i}=\vc{m}_0^i,\,\,\vc{m}|_{\Om_i}=\vc{m}^i,\,\,\vu_0|_{\Om_i}=\vu_0^i,\ \ \vu:=\f{\vc{m}}{R_0+Q_0},\,\,\vu^i:=\vu|_{\Om_i},$$ we conclude that there exist infinitely many weak solutions $[R_0,Q_0,\vc{m}]$ to the problem (\ref{intr1}) emanating from the initial data $[R_0,Q_0,\vu_0]$ and conserving the total energy. To see this, we need to verify the integral identities as follows.
\eqh{
&\intTO{\left(R_0 \p_t \phi +R_0 \vu \cdot \Grad_x \phi\right)}
=\int_0^T \!\!\!\int_{\cup_{i}\Om_i}\left(R_0 \p_t \phi +R_0 \vu \cdot \Grad_x \phi\right) \dxdt
\\
&=\sum_{i}\int_0^T \!\!\!\int_{\Om_i}\left(R_0^i \p_t \phi +R_0^i \vu^i \cdot \Grad_x \phi\right) \dxdt
=-\sum_{i} \int_{\Om_i}R_0^i \phi (0,\cdot)\dx
=-\int_{\Om} R_0 \phi(0,\cdot)\dx,
}
for any $\phi\in C_c^{\infty}( [0,T)\times \overline{\Om})$. This gives the weak formulation of the continuity equation (\ref{intr1})$_1$.

The  weak formulation of (\ref{intr1})$_2$ is verified in the same way. Let us now verify the momentum equation. It should be emphasized that we can choose $\chi$ to be the same sufficiently large positive constant in each subdomain $\Om_i$, in view of the assumptions of initial densities and the crucial observations (\ref{Inf8})-(\ref{Inf9}). Indeed, since $Z(R,Q)$ is uniformly bounded on $\Om$,
for any $i$ we can choose some $\chi$ large enough such that the right-hand side of (4.2)$_2$ is positive.
Therefore,
\eqh{
&\intTO{\Big( (R_0+Q_0) \vu\cdot \p_t \vc{\phi} +\left[(R_0+Q_0) \vu\otimes \vu + Z^{\gamma_{+}  } (R_0,Q_0) \mathbb{I}_3 \right]:\Grad_x \vc{\phi}\Big )}\\
&=\intTO{\left\{ (R_0+Q_0) \vu\cdot \p_t \vc{\phi} +\left[(R_0+Q_0)\vu\otimes \vu
+ Z^{\gamma_{+}  } (R_0,Q_0) \mathbb{I}_3
-\f{2}{3} \chi\,\mathbb{I}_3\right]:\Grad_x \vc{\phi}\right\}}
\\
&=\int_0^T \!\!\!\int_{\cup_{i}\Om_i}
\left\{ (R_0+Q_0) \vu\cdot \p_t \vc{\phi} +\left[(R_0+Q_0)\vu\otimes \vu
+ Z^{\gamma_{+}  } (R_0,Q_0) \mathbb{I}_3
-\f{2}{3} \chi\,\mathbb{I}_3\right]:\Grad_x \vc{\phi}\right\}
\dx\dt
\\
&=\sum_{i} \int_0^T \!\!\!\int_{\Om_i}
\left\{ (R_0^i+Q_0^i) \vu^i \cdot \p_t \vc{\phi} +\left[(R_0^i+Q_0^i )\vu^i \otimes \vu^i
+ Z^{\gamma_{+}  } (R_0^i,Q_0^i) \mathbb{I}_3
-\f{2}{3} \chi\,\mathbb{I}_3\right]:\Grad_x \vc{\phi}\right\}
\dx\dt
\\
&=-\sum_{i}\int_{\Om_i} (R_0^i +Q_0^i) \vu_0^i \cdot\vc{\phi}(0,\cdot)\dx
=-\int_{\Om}  (R_0+Q_0)   \vu_0 \cdot \vc{\phi}(0,\cdot)\dx,
}
for any $\vc{\phi}\in C_c^{\infty}( [0,T)\times \overline{\Om};\R^3)$ with $\vc{\phi} \cdot \vc{n}|_{\p \Om}=0$. Notice that $\vc{\phi} \cdot \vc{n}|_{\p \Om}=0$ allows us to write the first identity. In a similar manner, one gets the conservation of total energy by virtue of (\ref{lu91}). This finishes the proof of Theorem \ref{TH2}. $\Box$

\section{Further discussions} \label{fur}

\subsection{On other two-fluid models}\label{remo}
In this section, we advocate that Theorems \ref{TH1}-\ref{TH2} still hold for some other inviscid two-fluid models, by employing the same arguments as in Sections \ref{mawe}-\ref{maad}. The  governing equations still read as \eqref{intr1}, i.e.,
\begin{equation}\label{rem01}
\left\{\begin{aligned}
& \p_t R+\Div_x (R \vu)=0,\\
& \p_t Q+\Div_x (Q \vu)=0,\\
& \p_t [ (R+Q)\vu] +\Div_x [ (R+Q)\vu \otimes \vu ]+ \Grad_x p(R,Q)=\mathbf{0},
\end{aligned}\right.
\end{equation}
where $R$ and $Q$ are the densities of two fluids, and $\vu$ is the common velocity field, but the form of the pressure $p$ differs. We will choose two kinds of physically relevant models arising from engineering and asymptotic analysis respectively.

The first one is the liquid-gas flow, modeling the time evolution of the liquid and gas with the same velocity field. The model is widely applicable to describe the well and pipe flow processes \cite{EK1}. In this model, the pressure $P(m,n)$ is given by
\beq\label{remo2}
p(R,Q)=C \left(
-b(R,Q)+\sqrt{
b^2(R,Q)+c(R,Q)
}
\right),
\eeq
where
\begin{equation}\label{remo3}
\left\{\begin{aligned}
& b(R,Q)=k_0-R-a_0 Q, \\
& c(R,Q)=4k_0 a_0 Q,\\
\end{aligned}\right.
\end{equation}
and $C,k_0,a_0$ are positive physical constants. Mathematical results on the liquid-gas model are mainly concentrated on the viscous flows. Evje and Karlsen in \cite{EK1} proved the existence of global weak solutions to the one-dimensional viscous liquid-gas model by \emph{neglecting} the gas phase in the mixture momentum equation. The extension to the two-dimensional space was obtained by Yao et al. \cite{YZZ} with small initial data. In one-dimensional space, Evje et al. \cite{EWZ1} obtained the existence of global weak solutions by incorporating the gas phase in the mixture momentum equation. See also \cite{WYZ} and the references cited therein for more related results. On the other hand, again by \emph{neglecting} the gas phase in the mixture momentum equation, Ruan and Trakhinin in \cite{RuTr} introduced an entropy-like function and symmetrize the inviscid liquid-gas model, giving rise to many interesting results, say, the \emph{local-in-time} well-posedness of classical solutions. As far as we know, the existence of global-in-time solutions to the inviscid liquid-gas model (\ref{rem01}) in three-dimensional space is not known.

The second one is the fluid-particle model describing the time evolution of fluids and particles.   Under this setting, the pressure $p(R,Q)$ assumes the form
\beq\label{remo4}
p(R,Q)=R^{\gamma} + Q^{\beta}, \,\,\gamma, \beta \geq 1.
\eeq
Very recently, Vasseur et al. \cite{VWY} proved the existence of finite energy weak solutions to the viscous fluid-particle model under suitable constraints on $\gamma$ and $\beta$. For the inviscid fluid-particle model, we refer to Ruan and Trakhnin \cite{RuTr} for \emph{local-in-time} existence of shock waves and vortex sheets.

Similar to the two-fluid model (\ref{intr1}), we have the \emph{global-in-time} existence results for (\ref{rem01}) with pressure satisfying either (\ref{remo2})-(\ref{remo3}) or (\ref{remo4}). We leave the details to the interested reader.

\subsection{Local-in-time well-posedness} \label{local}
As we have seen in Theorem \ref{TH1}, problem (\ref{intr1})-(\ref{intr2}) is globally solvable but ill-posed in the class of weak solutions. However, it is locally well-posed when considering classical solutions. This is also in agreement with the compressible Euler equations. The main idea is to symmetrize
the inviscid two-fluid model (\ref{intr1}) and invoke the standard result of quasilinear symmetric hyperbolic system. Inspired by Ruan and Trakhinin \cite{RuTr} (see also Li et al. \cite{LSE1} for similar observation in 1D regime), we introduce the new variables. 
\[
\vr:= R+Q,\,\, s:= \f{R}{Q},
\]
and (\ref{intr1}) is then reformulated as
\begin{equation}\label{remo5}
\left\{\begin{aligned}
& \p_t \vr+\Div_x (\vr \vu)=0,\\
& \p_t s+\vu \cdot \Grad_x s=0,\\
& \p_t ( \vr \vu) +\Div_x ( \vr \vu \otimes \vu )+ \Grad_x p=\mathbf{0},\\
\end{aligned}\right.
\end{equation}
where the pressure $p$ is rewritten as 
\[
p=Z^{\gamma_{+}}(\vr,s).
\]
Here, $Z(\vr,s)$ is determined by
\begin{equation*}
\left\{\begin{aligned}
& \f{\vr}{1+s}=\left(1-\f{\vr s}{(1+s)Z} \right) Z^{\gamma}, \quad\gamma=\gamma_{+}/\gamma_{-}, \\
& \f{\vr s}{1+s} \leq Z,
\end{aligned}\right.
\end{equation*}
according to (\ref{intr2}). Suppose that $\vr$ and $s$ are strictly positive, one can then reformulate (\ref{remo5}) as 
\begin{equation}\label{remo6}
\left\{\begin{aligned}
& \f{1}{\vr \p_{\vr}p}(\p_t p+\vu \cdot \Grad_x p)+ \Div_x \vu=0,\\
& \p_t s+\vu \cdot \Grad_x s=0,\\
& \vr (\p_t \vu +\vu \cdot \Grad_x \vu)+ \Grad_x p=\mathbf{0}.\\
\end{aligned}\right.
\end{equation}
Notice that we may rewrite (\ref{remo6}) in the form of quasilinear symmetric hyperbolic system:
\beq\label{remo7}
\mathcal{A}_0(\mathbf{U}) \p_t \mathbf{U} +\sum_{i=1}^3 \mathcal{A}_i(\mathbf{U})
\p_{x_i} \mathbf{U}= \mathbf{0},
\eeq
where
\[
\mathbf{U}:=(p,\vu,s)^{t},
\]

\[
\mathcal{A}_0(\mathbf{U})=
\begin{pmatrix}
\f{1}{\vr \p_{\vr}p}&0&0&0&0\\
0&\vr &0&0&0 \\
0&0&\vr &0&0 \\
0&0&0&\vr  &0 \\
0&0 &0 &0&1 \\
\end{pmatrix},
\]

\[
\mathcal{A}_1(\mathbf{U})=
\begin{pmatrix}
\f{u_1}{\vr \p_{\vr}p}&1&0&0&0\\
1&\vr u_1&0&0&0 \\
0&0&\vr u_1 &0&0 \\
0&0&0&\vr u_1 &0 \\
0&0 &0 &0& u_1 \\
\end{pmatrix},
\]

\[
\mathcal{A}_2(\mathbf{U})=
\begin{pmatrix}
\f{u_2}{\vr \p_{\vr}p}&0&1&0&0\\
0&\vr u_2&0&0&0 \\
1&0&\vr u_2 &0&0 \\
0&0&0&\vr u_2 &0 \\
0&0 &0 &0& u_2 \\
\end{pmatrix},
\]

\[
\mathcal{A}_3(\mathbf{U})=
\begin{pmatrix}
\f{u_3}{\vr \p_{\vr}p}&0&0&1&0\\
0&\vr u_3&0&0&0 \\
0&0&\vr u_3 &0&0 \\
1&0&0&\vr u_3 &0 \\
0&0 &0 &0& u_3 \\
\end{pmatrix},
\]
Direct computation yields that
\[
\p_{\vr}p=
\frac{\gamma_{+}Z ^{\gamma_{+}}  (Z+s)  }
{\gamma (1+s) Z^{\gamma+1}  + \vr s},
\]
which is strictly positive if $\vr$ and $s$ are positive. This shows that $\mathcal{A}_0(\mathbf{U})$ is positive definite. Recalling the standard result on quasilinear symmetric hyperbolic system, see Kato \cite{Kato}, we conclude that the Cauchy problem for (\ref{remo7}) admits a unique local-in-time classical solution in Sobolev space $H^m(\R^3)$ for any integer $m>\f{5}{2}$.

\bigskip

\centerline{\bf Acknowledgement}
Y. Li is indebted to the Institute of Mathematics of the Czech Academy of Sciences for the invitation and hospitality from November 2018 to November 2019. 


\end{document}